\documentclass[12pt]{article}     \usepackage{amsmath}
\usepackage{amsthm}               \usepackage{amsopn}
\usepackage{psfig}                \usepackage{graphics}
\usepackage{latexsym}             \usepackage{amsfonts}
\usepackage{amssymb}
\swapnumbers
\theoremstyle{plain}
\newtheorem{theorem}{Theorem}[section]

\newtheorem{lemma}[theorem]{Lemma}
\newtheorem{corollary}[theorem]{Corollary}

\theoremstyle{definition}
\newtheorem{definition}[theorem]{Definition}
\newtheorem{remark}[theorem]{Remark}

\newenvironment{pf}{{\noindent\bf Proof.}}

\newcommand{\nc}[2]{ \newcommand{#1}{#2} }

\nc{\avint}{ {- \hspace{-3.5mm} \int} }  

\nc{\R}{{\rm {I \! R}}}  
\newcommand{\closure}[1]{ \stackrel{\rule{.1 in}{.01 in}}{#1} }

\newcommand{\chisub}[1]{ {\mathbf{\chi}}_{_{#1}} }

\newcommand{\newsec}[2]{ \section{#1} \label{sec-#2}  
                         \setcounter{equation}{0}     
                         \setcounter{theorem}{0} }    

\newcommand{\refeqn}[1]{ (\!\!~\ref{eq:#1}) } 
\newcommand{\refthm}[1]{ (\!\!~\ref{#1}) }    

\nc{\Holder}{H\"{o}lder\ }

\nc{\ith}{ \ensuremath{\text{i}^{\text{th}}} }
\nc{\jth}{ \ensuremath{\text{j}^{\text{th}}} }
\nc{\kth}{ \ensuremath{\text{k}^{\text{th}}} }
\nc{\curl}{ \nabla \times }
\nc{\Div}{ \nabla \cdot }

\begin{document}

\numberwithin{equation}{section}

\title{ The Hele-Shaw problem as a ``Mesa'' limit\\
        of Stefan problems: Existence, uniqueness,\\
        and regularity of the free boundary }
\author{Ivan A. Blank\\
        Marianne K. Korten\\
        Charles N. Moore\\
\normalsize}

\maketitle

\begin{abstract}
We study a Hele-Shaw problem with a mushy region obtained as a mesa
type limit of one phase Stefan problems in exterior domains.  We
study the convergence, determine some of the qualitative properties
and regularity of the unique limiting solution, and prove regularity
of the free boundary of this limit under very general conditions on
the initial data.  Indeed, our results handle changes in topology
and multiple injection slots.

\


\

\end{abstract}

\newsec{Introduction}{Intro}

Given a bounded domain $D \in \R^{n}$ with smooth boundary $\partial D,$ a 
finite set of closed curves $\{ s^{j}_{0}(x) = 0 \}$ such that $D$ is contained
in the union of their interiors, and a continuous function $p(x,t)$ defined on 
$\partial D \times (0, \infty),$ the classical Hele-Shaw problem with Dirichlet 
data is frequently formulated as follows:  Find a function $V(x,t)$ and a family 
of domains $S(t)$ (which each contain $D$) with $\partial S(t) = \{ t = s(x) \}$
and such that
\begin{equation}
\begin{array}{rll}
   \Delta_x V &= 0 \ \ &(x, \; t) \in S(t) \setminus D \\
   V(x, \; t) &= 0 \ \ &(x, \; t) \in \partial S(t) \\
   V(x, \; t) &= p(x) \ \ &(x, \; t) \in \partial D \times (0, \; \infty) \\
   \nabla_x V \cdot \nabla_x s &= \frac{\partial s}{\partial t} \ \ 
                               &(x, \; t) \in \partial S(t) \\
   \partial S(0) &= \cup_j \{ s^{j}_{0}(x) = 0 \} \;. \ \ & \ 
\end{array}
\label{eq:classicalHS}
\end{equation}

This problem models the advance of the slick formed by injecting oil between
two nearby plates, and has further been used in injection molding (used in
turn in the packaging industry, and more generally for the production of
plastic components, for example interior pieces of cars and aircraft), in
electrochemical machining (see \cite{MR}), and even to predict tumor growth
(see \cite{BF}).  Sometimes the normal derivative of $V$ at the ``slot,'' 
$\partial D,$ is prescribed, or curvature dependent terms are included in the
free boundary condition.  Among the most pressing open questions concerning
the Hele-Shaw problem are finding a weak formulation, studying the
regularity of $V(x,t)$ in $t,$ and determining the regularity of the free
boundary $t = s(x).$  Another question which has long attracted interest
is whether a Hele-Shaw problem could have a ``mushy'' region.

In turn, the ``Mesa'' problem describes the limit pattern 
$\lim_{m \rightarrow \infty} u_m$ of solutions $u_m$ of, say, the porous
medium equation, when the initial data are held fixed.  This problem
first appeared in connection with the modeling of problems related to
transistors (see \cite{EHKO}).  Caffarelli and Friedman studied the 
initial value problem in $\R^{n} \times [0,T]$ in \cite{CF}.  
They proved that the limit exists, that it is independent of the chosen 
subsequence, that it is independent of time, that it is equal to the 
characteristic function of one set plus the initial data times the 
characteristic function of the complement of that set, and finally
that that set can be characterized as the noncoincidence set of a 
variational inequality.  Further developments showed that the same 
conclusions hold for the limit when $u^m$ is replaced by a fairly 
general monotone constitutive function $\phi(u)$ with $\phi(0) = 0$ 
(see \cite{FH} for example), moreover, this behavior is a property of 
fairly general semigroups (see \cite{I}, \cite{BEG}).  A Mesa problem 
for an equation which gives rise to a mushy region was studied in 
\cite{BKM}.  For the Mesa problem we study in this paper we show all of 
the properties shown by \cite{CF} mentioned above, except that our 
limits will not be independent of time.  (This evolution in time is 
natural since we work in an outer domain where the inner boundary data 
serves as a source.)

In this paper the authors exploit a Mesa limit setting in an outer
domain $D^{c}$ to obtain naturally a weak formulation of the Hele-Shaw
problem (with Dirichlet condition as in Equation\refeqn{classicalHS}on 
the slot as above).  The use of one-phase Stefan problems with with 
``mushy'' regions and with increasing diffusivities naturally produces a mushy
region when we permit initial data $u_{I}$ for the approximating problems
to take values in the interval $[0, \; 1].$  These $u_{I}$ can be thought
of as generalized characteristic functions.  Another aspect of this approach
which is extremely attractive is the fact that changes in the topology of the 
``wet'' region do not interfere with the construction.  In short, whereas other
authors have taken a priori assumptions which ensure that their free boundary
stays smooth, we have been able to show existence of weak solutions for all
time, regardless of the possible changes in topology.  (Note that in \cite{CF} and
\cite{FH} they assume that their data is starlike with respect to the origin, and
in \cite{DL} log concavity of initial data is assumed to guarantee existence and
smoothness of the solutions.)  Indeed there are some very natural problems arising
in the applications where the topology should change.  Consider for example the 
problem of what happens with Hele-Shaw flow around an obstacle.  In this case, 
there is automatically a change in topology when the flow meets itself on the 
other side.  An interesting question is whether or not an air bubble will be left 
behind in the wake of the obstacle.  Another obvious problem from applications
where there will be changes in topology is if there are multiple injection slots.
In fact, this paper already deals with the second situation, since we never 
assume that $D$ is connected.

The classical version of the m-approximating problem which we use is given as follows:
\begin{definition}[m-approximating problem]   \label{mWholeMHS}
We let $u^{(m)}$ denote the solution of the following partial differential equation,
\begin{equation}
u^{(m)}_t = m \Delta (u^{(m)}(x,t) - 1)_{+}, \ \ \ \
            (x,t) \in D^{c} \times (0, \; + \infty) \cap \{ u^{(m)} > 1 \},
\label{eq:MesaHS}
\end{equation}
with boundary data given by
\begin{equation}
   \left.
   \begin{array}{rll}
      u^{(m)}(x,0) \! &= u_{I}(x), \; \ & x \in D^{c}, \\
      m(u^{(m)}(x,t) - 1)_{+} \! &= p(x), \; \
            & (x,t) \in \partial D \times (0, \; + \infty),
   \end{array}
   \ \ \ \ \right\}
\label{eq:bdryHS}
\end{equation}
and with the free boundary condition
\begin{equation}
    (-\nabla m(u^{(m)} - 1)_{+}, (1 - u_{I}) \;) \cdot \nu = 0 \;, \ \ 
    (x,t) \in \partial \{ u^{(m)} > 1 \} \;.
\label{eq:fbmHS}
\end{equation}
Here $\nu$ is the outer (n + 1) dimensional normal to the set $\{ u^{(m)} > 1 \},$
and this free boundary condition will be satisfied when the free boundary is smooth.
The weak formulation we give at the beginning of the next section will not require
any regularity assumptions on the free boundary or on the initial data.
For a fixed $m > 0,$ we call the free boundary problem determined by the equations
above the \textit{m-approximating problem}.  We will assume
\begin{equation}
   \begin{array}{l}
     0 \leq u_{I} \leq 1 \ \text{is compactly supported, } \\
     0 < p(x) \in C^{2, \alpha}(\partial D) \;, \\
     D \ \text{is a bounded set with} \ \partial D \in C^{2, \alpha}. \\
   \end{array}
\label{eq:makeass}
\end{equation}
\end{definition}
\noindent
\begin{remark}[Simplifying Assumptions]  \label{rmk:SA}
There are two assumptions which we make which serve to simplify the exposition but 
which are absolutely \textit{not} necessary for the derivation of our results:
\begin{enumerate}
   \item
The assumption that $u_{I}$ has compact support is only technical and serves to
simplify the construction of certain comparison functions in the fifth section.  We 
leave to the reader the verification that it suffices to assume only that 
$\{u_{I} = 1\}$ is compact for all of the results before the sixth section.
   \item
The assumption that the Dirichlet data, $p(x),$ on the boundary of the slot is a 
function of $x$ alone can be replaced with the assumption that the Dirichlet data
$p(x,t)$ is a nondecreasing function of $t$ for each fixed $x.$  On the other hand,
if $p(x,t)$ is allowed to decrease in time, then Lemma\refthm{Monoumint}does not 
hold in general, and this appears to be essential for the methods presented in this
paper.
\end{enumerate}
\end{remark}
For convenience, we will define
\begin{equation}
   M := ||p||_{L^{\infty}(\partial D)} \;.
\label{eq:Mdef}
\end{equation}
Now in terms of the positivity
assumption on $p(x),$ we note that $p(x) \equiv 0$ leads to a trivial case:
There will be no evolution at all.  To see this fact extend each $u^{(m)}$ by
$1$ across all of $D.$ In short, the positivity of $p(x)$ is driving the
evolution.

In this formulation $u^{(m)}$ is energy (or enthalpy), and $m(u^{(m)}(x,t) - 1)_{+}$
is temperature.  Equation\refeqn{MesaHS}comes from conservation of energy.
Basically, as $m$ increases, the diffusion happens faster.  Competing with this
increase in diffusion is the fact that the boundary data for $u^{(m)}$ on 
$\partial D$ is decreasing down to $1.$  As $m \rightarrow \infty$ we have
convergence of our operators to the following picture which is typical for Mesa
problems:
 
\

\psfig{file=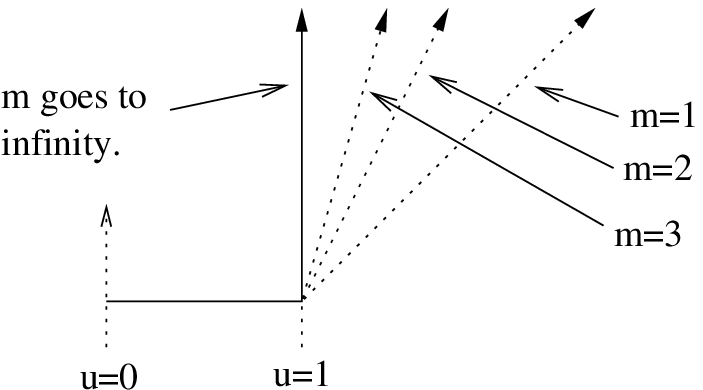}

\

We show that as $m \rightarrow \infty,$ the $u^{(m)}$ converge pointwise
to a limit $0 \leq u^{(\infty)} \leq 1,$ and $m(u^{(m)} - 1)_{+}$ tend 
pointwise to a bounded function $V(x, \; t).$  Furthermore, the function $V$ is
identically zero in $\{u^{(\infty)} < 1 \}$ and positive and harmonic for 
fixed $t$ in the component of the set $\{u^{(\infty)} = 1\}$ which contains $D.$ 
(On sets where $\{u^{(\infty)} = 1\}$ which are isolated from $D$ we will have
$V \equiv 0,$ see Remark\refthm{Vnc}\!\!.)
Finally, the pair $(u^{(\infty)}, \; V)$ solves our weak Hele-Shaw problem with 
a mushy region which we give formally in Definition\refthm{HSmr}\!\!.  In this way
we obtain the following results.
\begin{enumerate}
  \item We get a natural generalization of the weak formulation of the Hele-Shaw 
        problem of DiBenedetto and Friedman (see \cite{DF}).  Our formulation
        allows for a ``mushy'' region.  Moreover, after invoking a result
        of Bouillet, we will be able to say that solutions of our formulation
        are unique, and hence the Mesa limiting procedure we use gives an
        effective method of constructing the solution.
  \item We determine the regularity of the spatial slices of the free boundary 
        by invoking the results of Blank, Caffarelli, and Kinderlehrer and 
        Nirenberg for the regularity of the free boundary in the obstacle
        problem (see \cite{Bl}, \cite{C}, \cite{KN}).  Indeed, if we define
     \begin{equation}
        U(x,t) := \int_{0}^{t} u^{(\infty)}_{t}(x, \; s) \; ds
                = (1 - u_{I}(x))\chisub{A(t)}(x)
     \label{eq:Unowdef}
     \end{equation}
        where $A(t)$ gives the ``puddle'' at time $t,$ and
     \begin{equation}
        W(x,t) := \int_{0}^{t} V(x, \; s) \; ds \;,
     \label{eq:Wnowdef}
     \end{equation}
        and we formally apply the Baiocchi transformation to the equation
     \begin{equation}
        u^{(\infty)}_{t}(x, \; t) = \Delta_x V(x, \; t) \;,
     \label{eq:prelin}
     \end{equation}
        then for every fixed time $t_0,$ $W$ satisfies the obstacle problem:
     \begin{equation}
        0 \leq W(x,t_0) \;, \ \ \ \ \Delta_x W(x,t_0) = U(x,t_0) \;.
     \label{eq:UWsatOP}
     \end{equation}
\end{enumerate}
In a subsequent paper the authors will address the rectifiability and further space-time
regularity of the free boundary.

The paper is arranged as follows:  In section 2 we introduce our notion of weak
solutions to the approximating problem, show some qualitative properties of these
solutions, derive the existence of the limits (for now in weak-$* \ L^{\infty}$),
and give some trivial bounds on these limits.  Most of the results in this section
draw from the maximum principle and from the papers \cite{AK}, \cite{K1}, \cite{K2}, 
and \cite{K3}.  Section 3 gives a simple counter-example which motivates the definition 
of the free boundary and of the diffusive region.  In section 4 we derive some 
monotonicity properties of our sequences and limits, and as consequences we improve our 
weak-$* \ L^{\infty}$ convergence to pointwise convergence and give an explicit 
representation of $u^{(\infty)}.$  In section 5 we show that the limiting problem has a 
free boundary for all finite time (as opposed to if it ``escaped to infinity'' in zero 
time or by a fixed time), and we use explicit subsolutions to show that if an open ball 
within the diffusive region has a boundary point on the free boundary, then the free 
boundary will have nonzero velocity at that point.  In section 6 we apply the Baiocchi 
transformation to Equation\refeqn{prelin}and thereby derive a family of 
obstacle problems that this procedure yields.  At that point, under suitable assumptions
on the initial data, we can invoke the regularity theory for the obstacle problem (see 
\cite{Bl} and \cite{C}) to immediately derive regularity in space for the free boundary 
of our Hele-Shaw problem for almost every time.  In section 7, we extend the results of 
section 6 to include every time.  In the process of extending to every time, we establish
the continuity of the measure of the diffusive region with respect to time when the
initial data is strictly less than one.  Section 8 is an appendix which gives some of the 
barrier functions and some of their properties that we need in some of the earlier 
sections.

I. Blank wishes to thank Kansas State University for their hospitality while
the bulk of this paper was written.  M. Korten is indebted to Juan Luis 
V\'{a}zquez for introducing her to the Hele-Shaw problem.  M. Korten was partially 
supported by NSF EPSCoR Grant \#530517 under agreement NSF32169/KAN32170, 
and I. Blank's visits to Kansas State University were also supported by this 
grant.  All three authors wish to thank Barbara Korten for her patience and 
support while this paper was written.


\newsec{The weak formulations}{TWF}

For our m-approximating problem, we need an appropriate weak formulation which we
give here:

\begin{definition}[Weak solutions of the m-approximating problem] \label{WSmAP}
The nonnegative function $u^{(m)}(x, \; t) \in L^{1}_{loc}$ is a weak solution of
the m-approximating problem if for any
$\varphi \in C^{\infty}\left(\R^n \times [0, \infty)\right),$ such that
$\varphi \equiv 0$ on $\partial D \times [0, \infty),$ and $\varphi(x,t) \rightarrow 0$
as either $t \rightarrow \infty$ or $|x| \rightarrow \infty$ we have
\begin{equation}
\begin{array}{l}
 \displaystyle{\int_{D^c} \int_{0}^{\infty} \varphi_t(x,t) u^{(m)}(x,t) \; dt \; dx +
                 \int_{D^c} \int_{0}^{\infty} 
   \Delta_x \varphi(x,t) m \left[ u^{(m)}(x,t) - 1 \right]_{+} \; dt \; dx} \\
   \ \\
   = \displaystyle{\int_{\partial D} \int_{0}^{\infty}
        \frac{\partial \varphi}{\partial \nu}(x,t) p(x) \; dt \; d\mathcal{H}^{n - 1} x
                 -  \int_{D^c} \varphi(x,0) u_{I}(x) \; dx} \;.
\end{array}
\label{eq:wmap}
\end{equation}
\end{definition}  \noindent
We observe that the traces of $u^{(m)}$ and $m(u^{(m)} - 1)_{+}$ on the boundaries of
our domain will be well-defined by the work of Korten even if we only assume that we
are dealing with local solutions and that the initial trace is between zero and one.  
(See \cite{K1} which adapts the work of Dahlberg and Kenig, \cite{DK}, and see Lemma 
3.2 of \cite{K2}.)  Furthermore, we observe that the free boundary condition for the 
classical formulation (i.e. Equation\refeqn{fbmHS}\!) will be satisfied, whenever the 
functions and sets are sufficiently smooth.

For strictly local situations, it may be
helpful to note that if $\varphi \in C_{0}^{\infty}(\Omega)$ where
$\Omega \subset \subset D^c \times (0, + \infty),$ then we have
\begin{alignat*}{1}
   &\int_{D^c} \int_{0}^{\infty} \varphi_t(x,t) u^{(m)}(x,t) \; dt \; dx \\
   = - &\int_{D^c} \int_{0}^{\infty} m \Delta_x \varphi(x,t) \left[ u^{(m)}(x,t) - 1 \right]_{+}
                               \; dt \; dx \\
   = + &\int_{D^c} \int_{0}^{\infty} m \nabla_x \varphi(x,t) \cdot
                                     \nabla_x \left[ u^{(m)}(x,t) - 1 \right]_{+}
                               \; dt \; dx \;.
\end{alignat*}
The last integration by parts requires that we invoke the known regularity theory
and energy estimates for solutions of Equation\refeqn{MesaHS}\!\!.  See Lemma 1.2 of
\cite{AK}.

Since we want to discuss regularity of functions in Sobolev spaces, we need to fix ideas
about which representative we will use.  Although the ``Lebesgue point'' representative
where points are filled in by taking limits of averages over balls with radii going to
zero is the most common procedure, we will use a slightly different approach which exploits
the fact that the functions $u^{(m)}(x,t)$ are nondecreasing in $t$ for almost every $x$
by Lemma 4.2 of \cite{K1} and by the fact that the boundary data $p(x,t)$ is nondecreasing
in time for each fixed $x.$
\begin{definition}[Representative]   \label{fctrep}
For $t > 0,$ we set $$\tilde{u}^{(m)}(x, \; t) := \lim_{s \uparrow t} u^{(m)}(x, \; s) 
\ \ \text{for} \ a.e. \ x.$$
\end{definition}  \noindent
For economy of notation we will not bother to relabel any of our $u^{(m)}$'s, but always
understand that we are using $\tilde{u}^{(m)}$ as our representative for $u^{(m)}$ in all
pointwise matters.

We cannot expect $u^{(m)}$ to be continuous across the free boundary, and in
fact, the following theorem summarizes the results of Lemma 2.4 and Lemma 2.5 of \cite{K3}:
\begin{theorem}[$u^{(m)}$ must jump]  \label{ummj}
The set 
$$E := \{ x \in \R^n \; : 
       \; \exists \ t \ s.t. \; u_I(x) < u^{(m)}(x, \; t) < 1 \} $$
has Lebesgue $n$-dimensional measure equal to zero.
\end{theorem}
\noindent
On the other hand, $m(u^{(m)} - 1)_{+}$ will be continuous across the free boundary.  
(See \cite{K1} which verifies the required assumptions of \cite{DB} for continuity.)
For a fixed $m$ we now have the following picture of $u^{(m)}$ viewed from the side:

\

\psfig{file=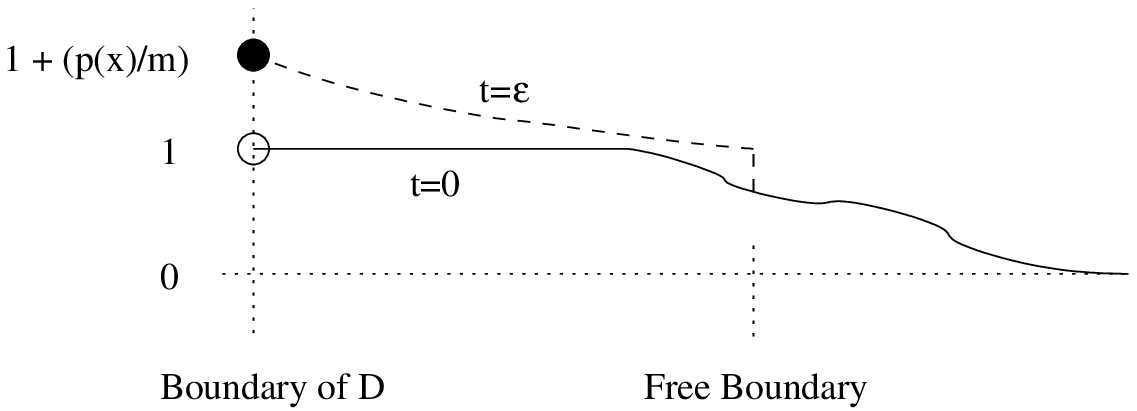}

\

By using the boundedness of the data along with the maximum principle that
the equation above enjoys, we can conclude that the solutions are bounded and
therefore by elementary functional analysis we can conclude that the following 
limit
\begin{equation}
         u^{(\infty)} := \lim_{m \rightarrow \infty} u^{(m)}(x,t), 
\label{eq:MesaLims}
\end{equation}
exists weak-$* \ L^{\infty}$ along a subsequence of $m$ in the entire domain.  
Again by using the maximum principle it is immediate that
\begin{equation}
         0 \leq u^{(\infty)} \leq 1 \;.
\label{eq:OneBdd}
\end{equation}
In fact, by using the maximum principle again, we can assert that
\begin{equation}
    0 \leq u^{(m)} \leq 1 + \frac{M}{m} \ \ \ \ \ \text{or} \ \ \ \ \ 
    0 \leq m(u^{(m)} - 1)_{+} \leq M \;,
\label{eq:mMfun}
\end{equation}
and we stress that $M$ is independent of $m$ and $t.$  Because of this fact, we
can take a further subsequence to ensure that $m(u^{(m)} - 1)_{+}$ has a 
limit $V$ in the weak-$* \ L^{\infty}$ topology of the entire domain.  

\begin{remark} \label{WooHooI} By combining Equation\refeqn{mMfun}with 
Theorem\refthm{ummj}we can conclude that the essential range of 
$u^{(m)}(x, \; \cdot)$ is a subset of $\{ u_I(x) \} \cup [1, \; 1 + M/m] \;.$
\end{remark}

Next, for $\varphi \in H^{1}_{0}(\Omega)$ we have:

\begin{alignat*}{1}
    \int_{0}^{+\infty} \int_{D^{c}} 
 \nabla_x \left[ m(u^{(m)} - 1)_{+} \right] \cdot \nabla_x \varphi \; dx \; dt
  & = \int_{0}^{+\infty} \int_{D^{c}} u^{(m)} \varphi_t \; dx \; dt \\
  & \rightarrow \int_{0}^{+\infty} \int_{D^{c}} u^{(\infty)} \varphi_t \; dx \; dt \;.
\end{alignat*}
Since $u^{(\infty)} \in L^{\infty}(\Omega) \subset L^{2}(\Omega),$ we have
\begin{equation}
   u^{(\infty)}_{t} \in H^{-1}(\Omega)
\label{eq:hmo}
\end{equation}
We summarize with the following lemma whose proof is now obtained trivially by
using the weak-$* \ L^{\infty}$ compactness we have and taking the limits as 
$m \rightarrow \infty$ in the weak formulation of the m-approximating problem. 
\begin{lemma}[Limits solve the limiting problem]  \label{limlimp}
Under our assumptions as above, we have 
\begin{equation}
   0 \leq V \leq M \;, 
\label{eq:VlimProb}
\end{equation}
and the pair $(u^{(\infty)}, \; V)$ satisfies
\begin{equation}
\!\!
\begin{array}{rl}
    &\displaystyle{\int_{D^c} \int_{0}^{\infty} \varphi_t(x,t) u^{(\infty)}(x,t) \; dt \; dx +
   \int_{D^c} \int_{0}^{\infty} \Delta_x \varphi(x,t) V(x,t) \; dt \; dx} \\
   \ \\
   = &\displaystyle{\int_{\partial D} \int_{0}^{\infty}
        \frac{\partial \varphi}{\partial \nu}(x,t) p(x) \; dt \; d\mathcal{H}^{n - 1} x 
    - \int_{D^c} \varphi(x,0) u_{I}(x) \; dx} \;.
\end{array}
\label{eq:whsprob}
\end{equation}
for any
$\varphi \in C^{\infty}\left(\R^n \times [0, \infty)\right),$ such that
$\varphi \equiv 0$ on $\partial D \times [0, \infty),$ and $\varphi(x,t) \rightarrow 0$
as either $t \rightarrow \infty$ or $|x| \rightarrow \infty.$
\end{lemma}
\begin{definition}[Hele-Shaw with mushy region]  \label{HSmr}
Any pair $(u^{(\infty)}, \; V)$ which satisfies Equation\refeqn{whsprob}for $\varphi$ as
in the lemma above will be called a weak solution of the Hele-Shaw problem with boundary
data $p(x)$ and initial data $u_{I}(x).$
\end{definition}
\begin{remark}[Inclusion of earlier models]  \label{Cwom}
The weak formulation we have above generalizes the Hele-Shaw formulation of DiBenedetto and
Friedman to include the situation where there is a mushy region (see \cite{DF}).  In
addition, our a priori regularity assumptions on the solutions are weaker.  (We only assume
$L^{1}_{loc}.$)
\end{remark}
\begin{theorem}[Uniqueness among all solutions]   \label{yitoo}
The solution of the limiting problem of Lemma\refthm{limlimp}is unique, so all
solutions of the Hele-Shaw problem as given in Definition\refthm{HSmr}are recoverable 
via the Mesa limit process we have introduced.
\end{theorem}
\begin{pf}
This theorem is an immediate application of the uniqueness theorem in section 3 
of \cite{Bo}.
\newline Q.E.D. \newline
\end{pf}


\newsec{Counter-example to regularity in time}{Cert}
Based on the assumptions we have so far, the function $V(x,t)$ will \textit{not} be
continuous in time in general.  To show this, we will assume for the sake of this
example that $u_{I}$ is a continuous function.  We will also assume in this section 
that the set $D$ is connected.  The trouble appears in certain cases where the set
\begin{equation}
   W := D \cup \{ x \in \R^n : u_{I}(x) = 1 \}
\label{eq:SetTroub}
\end{equation}
is disconnected.  If $K$ is a component of the set $\{ u_{I} = 1 \}$ and
$K$ is a positive distance away from $D,$ then the $u^{(m)}$ and therefore the
$u^{(\infty)}$ should not evolve on this set of $x$ until the free boundary
comes into contact with it.  Essentially, the ``patch'' $K$ will not ``see'' 
the input of the slot until the component of the set $\{ u^{(m)} \geq 1 \}$ 
which surrounds $D$ connects to it.
\begin{remark}[$V$ is not always continuous]  \label{Vnc}
To produce a situation where $V$ must be discontinuous in time,
simply consider the radially symmetric situation where $D = B_1$ and $u_{I}(r)$
is taken such that it is identically one on the set $3 \leq r \leq 5,$ but smaller
outside of it.  When what we want to call the free boundary reaches $r = 3,$ then 
it will instantaneously jump to $r = 5,$ and this leads to an immediate jump in 
the height of $V.$
\end{remark}
 
\
 
\psfig{file=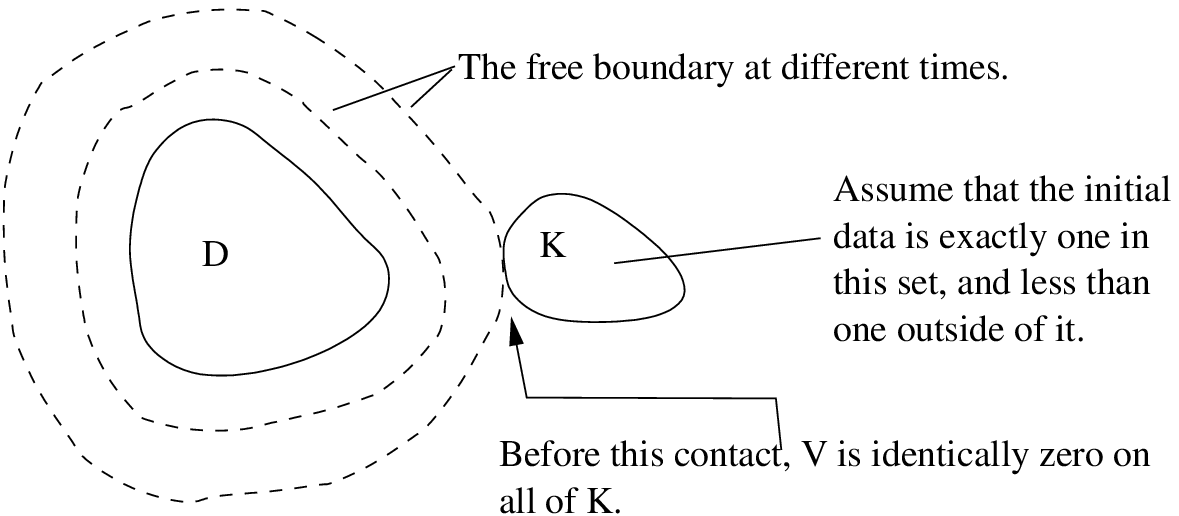}
 
\
 
\noindent
With these pictures in mind we make the following definitions:
\begin{definition}[Defining the free boundary]  \label{DFB}
We set
\begin{equation}
\begin{array}{l}
 t_{m}(x) := \inf \{ t \; : \; m(u^{(m)}(x,t) - 1)_{+} > 0 \} \\
 t_{\infty}(x) := \inf \{ t \; : \; V(x,t) > 0 \}
\end{array}
\label{eq:tdefs}
\end{equation}
and
\begin{equation}
\begin{array}{l}
 \tau_{m}(x) := \inf \{ t \; : \; u^{(m)}(x,t) \geq 1 \} \\
 \tau_{\infty}(x) := \inf \{ t \; : \; u^{(\infty)}(x,t) = 1 \} \;.
\end{array}
\label{eq:taudefs}
\end{equation}
We can also now define the diffusive region at time $t$ for the m-approximating
problem and the limiting Hele-Shaw problem to be respectively,
\begin{equation}
\begin{array}{rl}
  A^{(m)}(t) \!\! &:= \{ x \in D^c : \; m(u^{(m)}(x, \; t) - 1)_{+} > 0 \} \\
  A(t) \!\! &:= \{ x \in D^c \; : \; V(x, \; t) > 0 \} \;.
\end{array}
\label{eq:DiffReg}
\end{equation}
Now we simply define the free boundary at time $t$ to be 
$FB^{(m)}(t) := \partial A^{(m)}(t) \setminus \partial D,$ or 
$FB(t) := \partial A(t) \setminus \partial D,$ according to which problem we 
are considering.  (Here ``$\partial K$'' denotes the topological boundary of $K.$)
\end{definition}
\noindent
Observe that because of the radially symmetric example given above where $V$
is discontinuous in time, it is also clear that the free boundary cannot be
expected to vary continuously in time.
\begin{remark}[``Diffusive'' versus ``Wet'']    \label{DVW}
The reader may observe that ``Wet'' may be a more appropriate name for what we
are calling the diffusive region in the Hele-Shaw problem.  On the other hand,
``wet'' might more properly describe any region where $u^{(\infty)}(x,t) > 0,$
so we will stick to what is already appropriate in the approximating problems.
\end{remark}

By considering the figure above, it is clear that in general
$t_{m}(x)$ does not have to equal $\tau_{m}(x),$ and similarly with
$t_{\infty}(x)$ and $\tau_{\infty}(x).$  In the figure, they would
differ on the set $K.$  For $x \in K$ we have 
$\tau_{m}(x) = \tau_{\infty}(x) = 0,$ while $t_{m}(x)$ would be the
positive time when the moving part of the boundary of 
$\{u^{(m)} \geq 1 \}$ crosses $K,$ and $t_{\infty}(x)$ would be the
positive time when the moving part of the boundary of
$\{u^{(\infty)} = 1 \}$ crosses $K.$


\newsec{Monotonicity, inclusions, and consequences}{mic}
We start with a lemma which is a simple consequence of Lemma 4.2 in \cite{K1}.
\begin{lemma}[Monotonicity of $u^{(m)}$]   \label{Monoumint}
For each $m > 0$ and any $x \in D^c$ we have $u^{(m)}(x,t)$ is an increasing
function of time.
\end{lemma}
\begin{pf}
If $u^{(m)}(x,t) < 1,$ then this follows immediately from Lemma 4.2 of
\cite{K1}.  As a consequence of this fact, we can already say that the 
$A^{(m)}(t)$ are nested nondecreasing sets.  Now we consider the set 
$S := \{ (x,t) \; : \; u^{(m)}(x,t) > 1 \}.$  By observing that $u^{(m)}$
satisfies the heat equation within $S$ and invoking standard parabolic 
regularity theory, we can conclude that $u^{(m)}(x,t)$ will assume the
boundary values $1 + \frac{p(x)}{m}$ on $\partial D$ continuously, except 
possibly at the initial corner.  By Corollary 1.3  of \cite{AK} the function 
$(u^{(m)} - 1)_{+}$ is continuous within $D^c,$ so we can extend $u^{(m)}$ 
to be $1$ on the rest of the parabolic boundary of $S$ in a continuous 
fashion, except again at $t = 0$ where this set meets $\partial D.$  Now 
we consider the function 
\begin{equation}
    U(x,t) := u^{(m)}(x, \; t + \epsilon) - u^{(m)}(x, \; t) \;.
\label{eq:BUdef}
\end{equation}
Within $S,$ this function satisfies the heat equation.  On 
$\partial S \cap \partial D,$ we have $U(x,t) \equiv 0.$  (If we allow $p(x)$
to also depend on time, but we require it to be nondecreasing in time, then
we have $U(x,t) \geq 0$ here which is also fine.  It is not at all clear, 
however, how to deal with the most general case where $p(x,t)$ is allowed to
decrease in time.)  On the rest of the parabolic boundary of $S$ we have 
$U(x,t) = u^{(m)}(x, \; t + \epsilon) - 1 \geq 0.$  In the corner where
$U(x,t)$ is discontinuous, it stays bounded between zero and $1 + \frac{M}{m},$
and hence by the weak maximum principle, $U(x,t)$ is nonnegative in all of $D.$
Hence
\begin{equation}
    \frac{u^{(m)}(x, \; t + \epsilon) - u^{(m)}(x, \; t)}{\epsilon}
\label{eq:diffquotint}
\end{equation}
is nonnegative, and the result follows by taking $\epsilon \downarrow 0.$
\newline Q.E.D. \newline
\end{pf}

\begin{lemma}[Monotonicity of $u^{(\infty)}$]  \label{MonoUInf}
$u^{(\infty)}$ is monotone increasing in $t$ for $\mathcal{L}^n$ a.e. $x.$
\end{lemma}
\begin{pf}
By the last lemma the $u^{(m)}$ are monotone increasing in $t.$  Now to come to 
a contradiction, suppose $h > 0$ is a number and $\varphi(x,t)$ is a bounded 
nonnegative function which is compactly supported in our domain and which 
satisfies
  $$\int \int \varphi(x, \; t) u^{(\infty)}(x, \; t) \; dx \; dt >
    \int \int \varphi(x, \; t - h) u^{(\infty)}(x, \; t) \; dx \; dt \;.$$
Such a construction is possible if the conclusions of our theorem are not satisfied.
Then $\psi(x,t) := \varphi(x, \; t) - \varphi(x, \; t - h)$ is an admissible test
function which satisfies
  $$\int \int \psi(x, \; t) u^{(\infty)}(x, \; t) \; dx \; dt > 0 \;, \ \ \ \text{and} \ \ \
    \int \int \psi(x, \; t) u^{(m)}(x, \; t) \; dx \; dt \leq 0$$
which contradicts the weak-$* \ L^{\infty}$ convergence of the 
$u^{(m)}$ to $u^{(\infty)}.$
\newline Q.E.D. \newline
\end{pf}
Now without data on the slot which ``competes'' with the increasing diffusivity (or if
the functions solve the equation on all of $\R^n$), it is immediate by rescaling that
the diffusive regions (the $A^{(m)}(t)$) must be nondecreasing with $m.$  In our case, however, 
we need to produce an appropriate barrier.

\begin{theorem}[Temperature increases with $m.$]  \label{drim}
If $m < k,$ then
\begin{equation}
m[u^{(m)} - 1]_{+} \leq k[u^{(k)} - 1]_{+} \;.
\label{eq:tempinc}
\end{equation}
As a trivial consequence we can say that the $A^{(m)}(t)$ are nested:  
$m < k$ implies $A^{(m)}(t) \subset A^{(k)}(t) \;.$
\end{theorem}
\begin{pf}
We make the following definition:
\begin{equation}
   v^{(k)}(x, \; t) \; := \; \left\{
        \begin{array}{ll}
           u^{(m)}(x, \; t)  \ \ \ &\text{if} \ u^{(m)}(x, \; t) < 1 \\
           \ \\
           1 + \frac{m}{k}(u^{(m)}(x, \; t) - 1) \ \ \ &\text{if} \ u^{(m)}(x, \; t) \geq 1 \\
        \end{array}          \right.
\label{eq:sdvkdef}
\end{equation}
and claim that $v^{(k)}(x, \; t)$ is a subsolution of the k-approximating problem.
Indeed, this follows very quickly from the following two easily verifiable statements:
\begin{enumerate}
   \item The set where $v^{(k)} > 1$ is identical to the set where $u^{(m)} > 1.$
   \item $k(v^{(k)} - 1)_{+} \equiv m(u^{(m)} - 1)_{+}.$
\end{enumerate}
(Only the definition of $v^{(k)}$ is needed to verify these statements!)
By the first observation the free boundaries of $v^{(k)}$ and $u^{(m)}$ are identical 
in time and space, so the speed of these boundaries at every point is also identical.  
On the other hand, the second observation quickly leads to the conclusion that the function
$v^{(k)}$ satisfies the free boundary condition and the boundary value condition on the 
slot for the k-approximating problem exactly.  Finally, we simply compute:
\begin{alignat*}{1}
   v^{(k)}_{t}(x, \; t) & \leq u^{(m)}_{t}(x, \; t) \\
                        & = m \Delta (u^{(m)}(x, \; t) - 1)_{+} \\
                        & = k \Delta (v^{(k)}(x, \; t) - 1)_{+} 
\end{alignat*}
which shows that $v^{(k)}$ is a local subsolution to the k-approximating problem, and
therefore $v^{(k)}(x, \; t) \leq u^{(k)}(x, \; t).$  (In the first inequality above we 
have used Lemma\refthm{Monoumint}\!.)  We now have the simple consequence
\begin{equation}
   k[u^{(k)} - 1]_{+} \geq k[v^{(k)} - 1]_{+} \equiv m[u^{(m)} - 1]_{+} \;,
\label{eq:HGtm}
\end{equation}
which is the monotonicity we require.
\newline Q.E.D. \newline
\end{pf}
\begin{corollary}[Pointwise convergence of the temperature]  \label{pctemp}
The sequence of functions $\{m[u^{(m)} - 1]_{+}\}$ converges pointwise almost everywhere
to $V.$  In particular, the limiting function $V$ is unique.  (No subsequence is ever
needed.)
\end{corollary}
\begin{pf}
By the last theorem combined with the estimates in place already, at each point
we have a bounded increasing sequence.
\newline Q.E.D. \newline
\end{pf}
\begin{corollary}[Representation of $u^{(\infty)}$]  \label{uirep}
There is an increasing set-valued function of $t$ which we call $Q(t)$ such that
$u^{(\infty)}(x, \; t)$ admits the representation for almost every $(x, \; t)$:
\begin{equation}
   u^{(\infty)}(x,t) = \chisub{Q(t)}(x) + u_I(x) \chisub{Q(t)^{c}}(x) \;.
\label{eq:crurep}
\end{equation}
Furthermore, $u^{(\infty)}(x, \; t)$ is the pointwise limit of the functions
$u^{(m)}(x, \; t)$ almost everywhere, and $Q(t)$ can be chosen to be equal to 
the set $\{ x \in \R^{n} \; : \; \tau_{\infty}(x) < t \}$  ($Q(t)$ is increasing, means in 
the sense of set inclusion.)
\end{corollary}
\begin{pf}
By Lemma\refthm{MonoUInf}\!\!, we know that $u^{(\infty)}(x, \; t)$ is an increasing 
function of $t$ and so by using Lemma\refthm{MonoUInf}along with Equation\refeqn{OneBdd}we
can conclude that $u_{I}(x) \leq u^{(\infty)}(x, \; t) \leq 1.$  So, to prove 
Equation\refeqn{crurep}it suffices to show that $u^{(\infty)}(x, \; \cdot)$ does 
not attain values strictly between $u_{I}(x)$ and $1$ for a.e. $x.$  Indeed, by 
combining the previous theorem with Remark\refthm{WooHooI}we see that the limit
\begin{equation}
   \lim_{m \rightarrow \infty} u^{(m)}(x_0, \; t_0)
\label{eq:ptwiselim}
\end{equation}
exists for almost every $(x_0, t_0) \in \R^{n + 1}$ and is either equal to $1$ or to
$u_{I}(x_0).$  At that point it is a simple exercise to show that this pointwise limit
coincides with the weak-$* \ L^{\infty}$ limit almost everywhere.

In terms of the ``choice'' of $Q(t),$ it is clear that the only $x$ where there can 
be a choice is on the (possibly empty) set $\{ x \in \R^{n} \; : \; u_{I}(x) = 1 \}.$  
Examination of the definition of $Q(t)$ combined with the monotonicity of
$u^{(\infty)}(x, \; t)$ in $t$ makes it clear that 
$Q(t) := \{ x \in \R^{n} \; : \; \tau_{\infty}(x) < t \}$ will suffice.  (Note that
an equally good choice for $Q(t)$ would be to take 
$Q(t) := \{ x \in \R^{n} \; : \; t_{\infty}(x) < t \}.$)
\newline Q.E.D. \newline
\end{pf}
\begin{remark}[Uniqueness of limits of subsequences]  \label{Unique}
Although we had originally needed subsequences to be sure that our $u^{(m)}$ would
converge in weak-$* \ L^{\infty},$ the last result makes it clear that $u^{(m)}$ will
converge both pointwise and weak-$* \ L^{\infty}$ to $u^{(\infty)}$ without the need
to extract a subsequence.  
\end{remark}
\begin{corollary}[$V$ is harmonic within the diffusive region]  \label{DVEZ}
The spatial Laplacian of $V$ is zero for $x$ in the interior of $A(t).$  As a
consequence, if $\partial D \in C^{k, \alpha}$ where $k \geq 2,$ then $V$ attaches 
to the slot $D$ in a $C^{k, \alpha}$ fashion.
\end{corollary}
\noindent
Although the function $V$ is discontinuous in general, we do get regularity in space.
\begin{theorem}[Spatial continuity of $V.$]
$V(\; \cdot \;, \; t)$ is continuous for almost every time, $t.$
\end{theorem}
\begin{pf}
We will prove that $V(\; \cdot \;, \; t)$ is continuous by proving that it is both upper
and lower semicontinuous.  Lower semicontinuity follows immediately from the fact that
the functions $m[u^{(m)}(\; \cdot \;, \; t) - 1]_{+}$ are continuous and are increasing as 
functions of $m.$  To get upper semicontinuity we have a little bit more work.

By Lemma\refthm{MonoUInf}\!\!, $d\nu = u^{(\infty)}_t$ is a nonnegative measure.  Our Radon
measure can be ``sliced'' into Radon measures of one dimension less for a.e. t. (See
\cite{M} p. 139 - 142, and Equation (10.3) on p. 140 in particular.)  We will call
these slices $d\nu_t.$  Since for a.e. t we now have $\Delta_{x} V(x,t) = d\nu_{t},$
we know $V( \cdot, t)$ is subharmonic, and therefore upper semicontinuous for those
values of $t.$
\newline Q.E.D. \newline
\end{pf}


\newsec{Nontriviality and nondegeneracy}{NN}
Due to the competition between the decreasing data on $\partial D$ and the increasing
diffusivity of our m-approximating problems, we need now to rule out two trivial cases.
The possibilities that we need to exclude are:
\begin{enumerate}
   \item The possibility that the free boundary moves to infinity as soon as $t > 0.$
   \item The possibility that the free boundary never moves.
\end{enumerate}
Conceptually, our elimination of these cases is trivial.  We simply produce a family 
of supersolutions (and then subsolutions) to our m-approximating problems 
(Equation\refeqn{MesaHS}\!) whose free boundaries move in a suitable manner, 
independent of $m,$ and then we make use of the comparison principles that our
equations enjoy (see \cite{K1}).

\begin{theorem}[Boundedness of the free boundary]   \label{BFB}
If $u_{I}$ is compactly supported, then so is $u^{(\infty)}(x,t_0)$ for any $t_0 > 0.$
In fact, no matter what the initial data, there is always a supersolution whose free 
boundary has constant speed.
\end{theorem}
\noindent
Note that we need the compact support of the initial data for this result to hold.

\noindent
\begin{pf}
By using the boundedness of $D,$ and by translating and rescaling it will suffice to 
produce a supersolution to the problem where 
$D^c = B_1,$ \; $u_I = \chisub{ \{ B_2 \setminus B_1 \} },$ and $p(x) \equiv k.$
We claim that the function
\begin{equation}
   v^{(m)}(r,t) := \left( \frac{k(2 - r + \ell t)}{m(1 + \ell t)} + 1 \right)
                   \chisub{\displaystyle{(1, \; 2 + \ell t]}}(r)
\label{eq:vm}
\end{equation}
will suffice if $\ell$ is sufficiently large.  Note that this function is linear in the
radial variable, and that the free boundary moves with speed $\ell$ for all time.  The
following computations are elementary, and they prove our claim that $v^{(m)}$ is a
local supersolution:
\begin{equation}
  m \Delta (v^{(m)} - 1)_{+} = \frac{k(1 - n)}{r(1 + \ell t)} \leq 0 
                          \leq \frac{k \ell (r - 1)}{m(1 + \ell t)^2} = v^{(m)}_{t} \;.
\label{eq:locsuper}
\end{equation}
If we are on the free boundary, so that $r = 2 + \ell t,$ then
\begin{equation}
  \left| \frac{\partial}{\partial r} \left[ m(v^{(m)} - 1)_{+} \right] \right| = 
  \frac{k}{(1 + \ell t)} \leq \ell = \text{speed of the free boundary.}
\label{eq:fbsuper}
\end{equation}
\newline Q.E.D. \newline
\end{pf}
\begin{theorem}[Nontrivial motion of the free boundary]  \label{NMFB}
The boundary of the diffusive region does not remain stationary, and for any 
$R > 0,$ there exists a time $t_0$ such that for all $t > t_0,$ we have
$B_R \subset A(t_0).$
\end{theorem}
\begin{remark}   \label{EasWay}
If we only want to show that the free boundary does not remain stationary, then we can
use the same subsolutions as in the proof of Theorem\refthm{drim}\!.  On the other hand,
the subsolutions we use here are also needed in the proof of the next theorem.
\end{remark}
\begin{pf}
Using the fact that $D$ contains an open set, and by translating and rescaling it will 
suffice to produce a subsolution to the problem where (for any $\alpha \geq 0,$)
$D^c = B_{1 + \alpha},$ \; $u_I = \chisub{ \{ B_{2 + \alpha} \setminus B_{1 + \alpha} \} },$ 
and $p(x) \equiv k,$ and such that the free boundary of this subsolution moves a fixed
distance (independent of $\alpha$) in a finite time $T(\alpha).$

We will create our subsolution in a couple of steps.
First, let $w^{(m)}(r,t) = w^{(m)}_{1}(r,t) + w^{(m)}_{2}(r,t)$ where
$w^{(m)}_{1}(r,t)$ solves for each fixed $t$
\begin{equation}
  \begin{array}{rll}
    \Delta_x w^{(m)}_{1}(r,t) \! 
    &= 0 \ \ \ 
    & \text{in} \ B_{2 + \alpha + \ell t} \setminus B_{1 + \alpha} \\
    \ \\
    w^{(m)}_{1}(r,t) \! &= \displaystyle{1 + \frac{k}{m}} \ \ \ 
    & \text{on} \ \partial B_{1 + \alpha} \\
    \ \\
    w^{(m)}_{1}(r,t) \! &= 1 \ \ \ & \text{on} \ \partial B_{2 + \alpha + \ell t} \\
  \end{array}
\label{eq:wm1}
\end{equation}
and
$w^{(m)}_{2}(r,t)$ solves (again for each fixed $t$)
\begin{equation}
  \begin{array}{rll}
    \Delta_x w^{(m)}_{2}(r,t) \! 
    &= \displaystyle{\frac{\epsilon}{m}} \ \ \ 
    & \text{in} \ B_{2 + \alpha + \ell t} \setminus B_{1 + \alpha} \\
    \ \\
    w^{(m)}_{2}(r,t) \! 
    &= 0 \ \ \ 
    & \text{on} \ \partial \{ B_{2 + \alpha + \ell t} \setminus B_{1 + \alpha} \} \;.
  \end{array}
\label{eq:wm2}
\end{equation}
Since the free boundary will be given by $r = 2 + \alpha + \ell t,$ and we only need it
to move a fixed distance, we will also assume that $\ell t \leq 1.$  The functions 
$w^{(m)}_{1}$ and $w^{(m)}_{2}$ can be given explicitly, and their relevant properties
are given in the appendix.

By invoking Corollary\refthm{BOB}from the appendix, we have
\begin{equation}
  -\infty < -\tilde{C}_{1}(n)k \leq
  \frac{\partial}{\partial r} m w^{(m)}_{1}(2 + \alpha + \ell t, \; t) \leq
  -\tilde{C}_{2}(n)k < 0
\label{eq:Dmw1}
\end{equation}
where we stress that the constants are independent of $\alpha$ and $\ell t.$  
(See the appendix for a few more details, and remember that $0 \leq \ell t \leq 1.$)
By the same corollary we can conclude that
\begin{equation}
  0 < \tilde{C}_{3}(n) \epsilon \leq
  \frac{\partial}{\partial r} m w^{(m)}_{2}(2 + \alpha + \ell t, \; t) \leq
  \tilde{C}_{4}(n) \epsilon < \infty 
\label{eq:Dmw2}
\end{equation}
again with constants independent of $\alpha$ and $\ell t.$
Now we take $\epsilon > 0$ sufficiently small to ensure that
\begin{equation}
  -\infty < -C_1 k \leq
  \frac{\partial}{\partial r} \left[ m w^{(m)}(2 + \alpha + \ell t, \; t) \right] \leq
  -C_2 k < 0 \;.
\label{eq:Dmw}
\end{equation}
Now by taking $\ell \leq C_2 k/2$ we can be sure that our function is a subsolution
along the free boundary.  We note that 
$\Delta m(w^{(m)}(r, \; t) - 1)_{+} = \epsilon > 0$ in the region where
$w^{(m)}(r, \; t) > 1.$  Since
\begin{equation}
   \lim_{m \rightarrow \infty} w^{(m)}_t(r, \; t) = 0 \;,
\label{eq:wtmLIM}
\end{equation}
Once $m$ is sufficiently large, we automatically have
\begin{equation}
   \Delta m (w^{(m)}(r, \; t) - 1)_{+} \geq w^{(m)}_t(r, \; t) \;.
\label{eq:wmlocsub}
\end{equation}
Since the free boundary of our subsolution moves with speed $\ell > 0,$ we are done.
\newline Q.E.D. \newline
\end{pf}

The following theorem shows the instantaneous detachment of the free boundary
from the slot, $\partial D,$ even if $u_{I}(x) \leq \lambda < 1$ in all of $D^{c}.$
\begin{theorem}[Instantaneous formation of the diffusive region]  \label{IFDR}
If $t_0 > 0,$ The set $A(t_0)$ contains an open neighborhood of $\partial D.$
\end{theorem}
\begin{pf}
Because we have assumed that $\partial D \in C^{2, \alpha} ,$ every point on 
$\partial D$ can be touched from within with a tangent ball, and then we can use
the same subsolutions of the previous theorem to force instantaneous movement of
the free boundary.
\newline Q.E.D. \newline
\end{pf}


\newsec{Spatial regularity results}{SRFB}
In this section we will derive spatial regularity for both the function
$V(x,t)$ and for the free boundary.  We apply the Baiocchi 
transformation to $V(x,t)$ and define:
\begin{equation}
     W(x,t) := \int_{0}^{t} V(x, \; s) \; ds \;.
\label{eq:WBaio}
\end{equation}
Observe that by Lemma\refthm{MonoUInf}(which shows that the diffusive regions
are increasing in time) and by the positivity of $V(x,t)$ in the diffusive
region, it is clear that the set $\{ W > 0 \}$ is identical to the set
$\{ V > 0 \}.$  Now to find the regularity of 
$\partial \{ x \in \R^{n} \; : \; W(x, T) > 0 \} $ we will show that 
$W(\cdot,t)$ belongs to $H^1_{loc}(D^{c})$ for almost every $t,$ and 
then that $W(x, t)$ is a weak solution of the following obstacle 
problem in almost every time slice $t = T$
\begin{equation}
   0 \leq W(x, T), \ \ \ \ \ 
   \Delta_x W(x, T) = \chisub{ \{ W(x, T) > 0 \} }(x)(1 - u_{I}(x)) \;.
\label{eq:WobProb}
\end{equation}
After that we will be able to invoke regularity results for the obstacle problem 
due to Caffarelli, Kinderlehrer, Nirenberg, and Blank.

For simplicity, we let $\alpha_{m}(s) := m(s - 1)_{+}.$  We start by stating some
simple trace results.  Basically, we need to adapt Equation\refeqn{wmap}to some
situations with slightly different test functions.

\begin{lemma}[First Trace result] \label{TrRes}
If $\psi(x) \in C^{\infty}$ is supported in the interior of $D^{c},$ then the following formula
holds for a.e. $T$:
\begin{equation}
\begin{array}{l}
 \displaystyle{\int_{D^c} \int_{0}^{T}
   [\Delta_x \psi(x)] \; \alpha_{m}(u^{(m)}(x,t)) \; dt \; dx} \\
   \ \\
   = \displaystyle{ \int_{D^c} \psi(x) [u^{(m)}(x,T) - u_{I}(x)] \; dx} \;.
\end{array}
\label{eq:wmap2}
\end{equation}
\end{lemma}
\begin{pf}
We make the following definition
\begin{equation}
   \varphi(x,t) := \left\{
       \begin{array}{rl}
          \psi(x) \ \ \ & t \leq T \\
          \ \\
          \Theta(x,t) \ \ \ & t > T \\
       \end{array}
   \right.
\label{eq:nicetest}
\end{equation}
where $\Theta(x,t)$ is chosen to ensure that $\varphi(x,t)$ is a permissible
test function for our m-approximating problem.  In particular, we need
$\varphi(x,t) \in C^{\infty}$ and we need it to converge to zero as
$t \rightarrow \infty.$  Neither requirement poses any difficulty.

By using the trace result of \cite{AK} (see Theorem 1.1 of \cite{AK}) our
functions $u^{(m)}(x,t)$ solve our m-approximating problem starting at time $T$
with initial data $u^{(m)}(x,T)$ for almost every $T,$ and so (for those $T$)
we can use $\Theta(x,t)$ as the test function in Equation\refeqn{wmap}to obtain:
\begin{equation}
\begin{array}{rl}
 \ & \displaystyle{\int_{D^c} \int_{T}^{\infty} \Theta_t(x,t) u^{(m)}(x,t) \; dt \; dx} \\
   \ \\
 + & \displaystyle{\int_{D^c} \int_{T}^{\infty}
   \Delta_x \Theta(x,t) \alpha_{m}(u^{(m)}(x,t)) \; dt \; dx} \\
   \ \\
 = - & \displaystyle{\int_{D^c} \psi(x) u^{(m)}(x,T) \; dx} \;.
\end{array}
\label{eq:wmapT}
\end{equation}
(Note that $\Theta(x,T) = \psi(x)$ since we have required that $\varphi(x,t)$ be
smooth.)  Now by subtracting this equation from what we have when we plug in the
function $\varphi(x,t)$ defined in Equation\refeqn{nicetest}into
Equation\refeqn{wmap}we get Equation\refeqn{wmap2}immediately.
\newline Q.E.D. \newline
\end{pf}

The proof of the following result is almost identical to the proof of the trace result
above, so we omit it.
\begin{lemma}[Second Trace result] \label{TrRes2}
For a.e. $t_0, t_1$ such that $0 \leq t_0 < t_1 < \infty,$ and for 
$\varphi \in C^{\infty}\left(\R^n \times [t_0,t_1]\right),$ which satisfies
$\varphi \equiv 0$ on an open set containing $D \times [t_0,t_1],$ we have
\begin{equation}
\begin{array}{l}
 \displaystyle{\int_{D^c} \int_{t_0}^{t_1} \varphi_t(x,t) u^{(m)}(x,t) \; dt \; dx +
                 \int_{D^c} \int_{t_0}^{t_1}
   \Delta_x \varphi(x,t) \alpha_{m}(u^{(m)}(x,t)) \; dt \; dx} \\
   \ \\
   = \displaystyle{\int_{D^c} \left( \left. \rule[-.02 in]{.0 in}{.22 in}
               \left[  \varphi(x,s) u^{(m)}(x,s) \right] \right|_{s = t_0}^{s = t_1} \right) \; dx} \;.
\end{array}
\label{eq:wmap3}
\end{equation}
\end{lemma}

Now we state the standard energy estimate for our situation.
\begin{lemma}[Energy estimates for the $u^{(m)}$]   \label{umenergy}
Let $0 < r < R$ and $0 \leq t_0 < t_1.$  Then there is a constant of
the form
\begin{equation}
   C = \frac{C(n)}{(R - r)^2}
\label{eq:rtdep}
\end{equation}
such that the following energy estimate holds:
\begin{equation}
    \int_{t_0}^{t_1} \int_{B_{r}(x_0)} 
    | \nabla \alpha_{m}(u^{(m)}) |^2 \; dx \; dt \leq
    C\int_{t_0}^{t_1} \int_{B_{R}(x_0)} \alpha_{m}(u^{(m)})^2 \; dx \; dt \;.
\label{eq:umspace}
\end{equation}
\end{lemma}
\begin{remark}[Independence of Time]   \label{rmk:IoT}
Notice that the constant $C$ is independent of time and notice that the time
intervals in the integrals in each side of the inequality are identical.
\end{remark}
\begin{pf}
We choose $\eta(x) \in C^{\infty}_{0}(B_R)$ such that
\begin{enumerate}
   \item $\eta \equiv 1$ on $\closure{B_{r}},$
   \item $0 \leq \eta \leq 1,$ and
   \item $|\nabla \eta| \leq 4(R - r)^{-1} \;.$
\end{enumerate}
Now let $\varphi(x,t) := \alpha_{m}(u^{(m)}) \eta(x)^2,$ and apply the last lemma
and Green's identity to obtain:
\begin{equation}
\begin{array}{l}
 \displaystyle{\int_{B_R} \int_{t_0}^{t_1} \varphi_t(x,t) u^{(m)}(x,t) \; dt \; dx -
               \int_{B_R} \left( \left. \rule[-.02 in]{.0 in}{.22 in} 
               \left[  \varphi(x,s) u^{(m)}(x,s) \right] \right|_{s = t_0}^{s = t_1} \right) \; dx} \\
   \ \\
   = \displaystyle{\int_{B_R} \int_{t_0}^{t_1} \nabla_{x} \varphi(x,t) \nabla_{x} \alpha_{m}(u^{(m)}) \; dt \; dx} \;.
\end{array}
\label{eq:prehappy}
\end{equation}
The fact that we can apply Green's identity above can be justified using Lemma 1.2
of \cite{AK}.  By Lemma\refthm{Monoumint} we know that for a.e. $x \in D^c, \ \mu_x(t):=u_t(x,t)$
is a non-negative Radon measure, whence it suffices to have $\phi(t) \in C_0(0, +\infty)$ for the
distributional pairing $(\mu_x(t), \phi(t))$ to be defined, and non-negative if $\phi \geq 0.$
Using Equation\refeqn{wmap}it is easy to see that the function $g(x) := (\mu_x(t), \chisub{[t_0,t_1]}(t))$
is locally integrable in $x.$  In other words, we can integrate the left hand side of
Equation\refeqn{prehappy}by parts in time to give us the following inequality:
\begin{alignat*}{1}
  \ & \ \ \ \int_{B_R} \int_{t_0}^{t_1} \nabla_{x} \left( \alpha_{m}(u^{(m)}) \eta^2(x) \right) 
                                    \nabla_{x} \alpha_{m}(u^{(m)}) \; dt \; dx \\
  \ &= \; - \int_{B_R} \left[ \int_{t_0}^{t_1} \left( \alpha_{m}(u^{(m)}) \eta^2(x) \right) 
                                           d\mu_x(t) \right] \; dx \\
  \ &\leq 0.
\end{alignat*}
This inequality implies
$$\begin{array}{c}
    \displaystyle{
\int_{B_R} \int_{t_0}^{t_1} \eta^2 \left| \nabla_{x} \alpha_{m}(u^{(m)}) \right|^2 \; dt \; dx \leq} \\
    \ \\
    \displaystyle{
  2\int_{B_R} \int_{t_0}^{t_1} \left( \eta \left| \nabla_{x} \alpha_{m}(u^{(m)}) \right| \right)
                              \left( \alpha_{m}(u^{(m)}) \left| \nabla_{x} \eta \right| \right)
                                    \; dt \; dx}
  \end{array}$$
and so Cauchy-Schwarz gives
\begin{equation}
\begin{array}{l}
 \displaystyle{\int_{t_0}^{t_1} \left[ 
    \int_{B_R} \eta(x)^2 \left| \nabla_{x} \alpha_{m}(u^{(m)}) \right|^2 dx \; \right] dt } \\
   \ \\
   \leq \displaystyle{4\int_{t_0}^{t_1} \left[
    \int_{B_R} \left| \nabla_{x} \eta(x) \right|^2 \alpha_{m}(u^{(m)})^2 dx \; \right] dt } \;.
\end{array}
\label{eq:extrahappy}
\end{equation}
\newline Q.E.D. \newline
\end{pf}
The aforementioned independence of time of the constant $C$ leads immediately to the following
theorem:
\begin{theorem}[Time-slice energy estimates] \label{T-seeum}
For almost every time $t$ we have
\begin{equation}
    \int_{B_{r}(x_0)} | \nabla \alpha_{m}(u^{(m)}(x,t)) |^2 \; dx \leq
    C\int_{B_{R}(x_0)} \alpha_{m}(u^{(m)}(x,t))^2 \; dx \;,
\label{eq:umspace2}
\end{equation}
and also
\begin{equation}
    \int_{B_{r}(x_0)} | \nabla V(x,t) |^2 \; dx \leq
    C\int_{B_{R}(x_0)} V(x,t)^2 \; dx \;,
\label{eq:VspaceCacc}
\end{equation}
where once again, the constant $C$ has the form given in Equation\refeqn{rtdep}\!.
\end{theorem}
\begin{pf}
Equation\refeqn{umspace2}is already known from the previous lemma, and Equation\refeqn{VspaceCacc}follows
from Equation\refeqn{umspace2}by using the $L^2$ convergence of $\alpha_{m}(u^{(m)}(x,t))$ to
$V(x,t)$ and by using the lower semicontinuity of the Dirichlet integral.
\newline Q.E.D. \newline
\end{pf}
\begin{corollary}[Energy estimates for $V(\cdot,t)$]  \label{VLip}
For almost every $t, \ V(\cdot,t) \in H^{1}_{loc}(D^c).$
\end{corollary}
\begin{theorem}[Energy estimates for $W(\cdot,t)$]  \label{WH1inSpace}
For all $t > 0, \ W(\cdot,t) \in H^{1}_{loc}(D^c).$
\end{theorem}
\begin{pf}
Let $K$ be a compact subset of $D^{c}.$  Since $W(\cdot,t) \leq Mt,$ we have
$W(\cdot,t) \in L^{\infty}(K) \subset L^{2}(K).$  It remains to show that
$\frac{\partial W}{\partial x_i}(\cdot,t) \in L^{2}(K).$

Let $\varphi \in C_{0}^{1}(K)$  
and estimate.
  \begin{alignat*}{1}
  \left( \frac{\partial W}{\partial x_i}, \; \varphi \right) \;
  &= \; - \int_{K} \; W(x,t) \frac{\partial \varphi}{\partial x_i}(x) \; dx \\
  &= \; - \int_{K} \int_{0}^{t} V(x, \; s) \frac{\partial \varphi}{\partial x_i}(x) \; ds \; dx \\
  &= \; - \int_{K} \int_{0}^{t} \lim_{m \rightarrow +\infty} \alpha_{m}(u^{(m)}(x, \; s)) \;
          \frac{\partial \varphi}{\partial x_i}(x) \; ds \; dx \\
  &= \; - \lim_{m \rightarrow +\infty} \int_{K} \int_{0}^{t} \alpha_{m}(u^{(m)}(x, \; s)) \;
          \frac{\partial \varphi}{\partial x_i}(x) \; ds \; dx \\
  &= \; \lim_{m \rightarrow +\infty} \int_{K} \int_{0}^{t} 
        \frac{\partial}{\partial x_i} \left( \alpha_{m}(u^{(m)}(x, \; s)) \right) \varphi(x) \; ds \; dx \\
  &= \; \lim_{m \rightarrow +\infty} \int_{K} \varphi(x) \left[ \int_{0}^{t}
        \frac{\partial}{\partial x_i} \left( \alpha_{m}(u^{(m)}(x, \; s)) \right) \; ds \right] \; dx \\
  \end{alignat*}
where we have used Lebesgue's Dominated Convergence Theorem repeatedly in the computation above.
Now, by using Minkowski's integral inequality and the Cauchy-Schwarz Inequality we have
  \begin{alignat*}{1}
  \ & \left| \int_{K} \varphi(x) \left[ \int_{0}^{t} \frac{\partial}{\partial x_i} 
      \left( \alpha_{m}(u^{(m)}(x, \; s)) \right) \; ds \right] \; dx \right| \; \\
  &\leq \left| \int_{K} \varphi(x)^2 \; dx \right|^{1/2} \left| \int_{K} \left(
        \int_{0}^{t} \frac{\partial}{\partial x_i} 
        \left( \alpha_{m}(u^{(m)}(x, \; s)) \right) \; ds \right)^{2} dx \right|^{1/2} \\
  &\leq ||\varphi||_{L^2(K)} \int_{0}^{t} \left| \int_{K} \left( \frac{\partial}{\partial x_i}
        \left( \alpha_{m}(u^{(m)}(x, \; s)) \right) \right)^2 dx \right|^{1/2} ds \\
  &\leq ||\varphi||_{L^2(K)} \left| \int_{0}^{t} \int_{K} \left( \frac{\partial}{\partial x_i}
        \left( \alpha_{m}(u^{(m)}(x, \; s)) \right) \right)^2 dx \; ds \right|^{1/2} \cdot
        \left[ \int_{0}^{t} ds \right]^{1/2} \\
  &\leq t^{1/2} ||\varphi||_{L^2(K)} ||\nabla_{x} \alpha_{m}(u^{(m)}(x, \; s)) ||_{L^2(K \times (0,t))} \\
  &\leq C t^{1/2} ||\varphi||_{L^2(K)}
  \end{alignat*}
where the last constant is independent of $m$ by Lemma\refthm{umenergy}\!.
\newline Q.E.D. \newline
\end{pf}
\begin{remark}   \label{WC1al}
In fact, once we show that $W(\cdot,t)$ is a solution of the obstacle problem, we
will be able to infer from elliptic regularity theory, that $W(\cdot,t)$ is 
$C^{1,\alpha}$ in space for all $\alpha < 1.$
\end{remark}
In order to derive Equation\refeqn{WobProb}we will need to commute the 
Laplacian with the integral in time.  To accomplish this commutation we turn 
back to the approximating problem where this switch is simpler. 

By using weak-$* \ L^{\infty}$ convergence of the temperature functions we have
\begin{alignat*}{1}
   \lim_{m \rightarrow \infty} &\int_{D^c} \int_{0}^{T}
    [\Delta_x \psi(x)] \; m \left[ u^{(m)}(x,t) - 1 \right]_{+} \; dt \; dx \\
   = &\int_{D^c} \int_{0}^{T} [\Delta_x \psi(x)] \; V(x, t) \; dt \; dx \\
   = &\int_{D^c} [\Delta_x \psi(x)] W(x,T) \; dx \;.
\end{alignat*}
On the other hand, by using the dominated convergence theorem we conclude that
for almost every $T$
\begin{alignat*}{1}
   \lim_{m \rightarrow \infty} &\int_{D^c} \psi(x) [u^{(m)}(x,T) - u_{I}(x)] \; dx \\
   = &\int_{D^c} \psi(x) [u^{(\infty)}(x,T) - u_{I}(x)] \; dx \\
   = &\int_{D^c} \psi(x) \chisub{ \{ W(x, T) > 0 \} }(x)(1 - u_{I}(x))
\end{alignat*}

Since $W \geq 0,$ we can combine the last two computations with the previous
lemma to conclude that $W(\cdot, \; t)$ solves the obstacle problem for a.e. $t > 0.$
Equation\refeqn{WobProb}and this fact allows us to use the technology of 
\cite{Bl} and \cite{C} to infer regularity of the free boundary here as long 
as we satisfy the nondegeneracy condition:
\begin{equation}
   u_I(x) \leq \lambda < 1 \;.
\label{eq:ndsit}
\end{equation}
($u_I$ satisfying Equation\refeqn{ndsit}will be referred to as \textit{nondegenerate
initial data}.)  In this case we will have the next theorem which we state after one 
simple definition.
\begin{definition}[Minimum diameter]  \label{mindiam}
The minimum diameter of a set $S \subset \R^n$ (denoted ``$m.d.(S)$'') is the
infimum of the distances between parallel hyperplanes enclosing $S.$
\end{definition}
\begin{theorem}[Regularity of the free boundary in space]  \label{RFBS}
Assume $u_{I}$ is continuous, nondegenerate initial data.  Then there is a modulus 
of continuity $\sigma$ which depends on $\lambda, n,$ and the modulus of continuity 
of $u_{I}$ such that at almost any time $t,$ (indeed every time $t_0$ where 
Equation\refeqn{WobProb}is valid) and for any free boundary point $(x_0,t_0)$ 
contained in the interior of $D^{c}$ we have either (with $B_r(x)$ denoting the 
(spatial) ball centered at $x$ with radius $r$)
\begin{equation}
   m.d.(B_r(x_0) \cap A(t)^{c}) \leq r \sigma(r) \ \ \ \ \text{for all} \ r \leq 1 \;,
\label{eq:singularpt}
\end{equation}
(so that the nondiffusive region is ``cusp-like'') or
\begin{equation}
   \lim_{r \rightarrow 0} \frac{|B_r(x_0) \cap A(t)|}{|B_r|} = \frac{1}{2} 
\label{eq:regpt}
\end{equation}
in which case we will say that $(x_0, t_0)$ is a regular point of the free boundary.
\end{theorem}
\begin{pf}
Simply use Equation\refeqn{WobProb}together with the results from \cite{Bl}.
\newline Q.E.D. \newline
\end{pf}
\begin{theorem}[Better regularity] \label{betreg}
Assume the hypotheses of the previous theorem and assume that $(x_0, t_0)$ is a
regular point of the free boundary.  Then the free boundary intersected with
$\{t = t_0\}$ will be Reifenberg vanishing near $x_0.$  (For a definition of
Reifenberg vanishing, see \cite{Bl}.)  If $u_I$ is Dini continuous, then the
free boundary will be $C^1$ (in space) near $x_0,$ and if $u_I \in C^{k,\alpha}$ 
then the free boundary will be $C^{k+1,\alpha}$ (in space) near $x_0.$
\end{theorem}
\begin{pf}
Again, just combine Equation\refeqn{WobProb}with the results from \cite{Bl}.
\newline Q.E.D. \newline
\end{pf}
\begin{remark}[Clarification]  \label{clear}
When we say ``near'' $x_0$ in the theorem above we mean near in space only.
In other words, we are specifically talking about the free boundary restricted
to the time slice $t = t_0.$
\end{remark}
At this point we need to call attention to our assumptions that $u_I$ is continuous
and nondegenerate.  Indeed it is well known that there are examples of ``persistent 
corners'' in Hele-Shaw problems.  (See \cite{KLV} for example.)  Our theorem does 
not contradict this fact.  Our spatial regularity theorem says nothing if $u_I$ is 
discontinuous, or if the set $\{ u_I = 1 \} \cap D^{c}$ is nonempty.  On the other 
hand, our assumptions do not rule out ``focusing'' or changes in topology of the 
diffusive region, so it is certainly nontrivial that corners do not arise in this 
setting.  Consider a case where $u_I$ is very close to zero except for an $\epsilon$ 
neighborhood of an annulus which contains the slot.  In this annulus assume that 
$u_I$ is extremely close to one.  (The epsilon neighborhood is needed to make $u_I$
continuous.)  In this case, the diffusive region will make its way around the annulus
in each direction very quickly, but expand slowly into the region inside and outside
of the annulus.  Thus, it will meet itself on the other side of the annulus long before
it fills in the interior.

\

\psfig{file=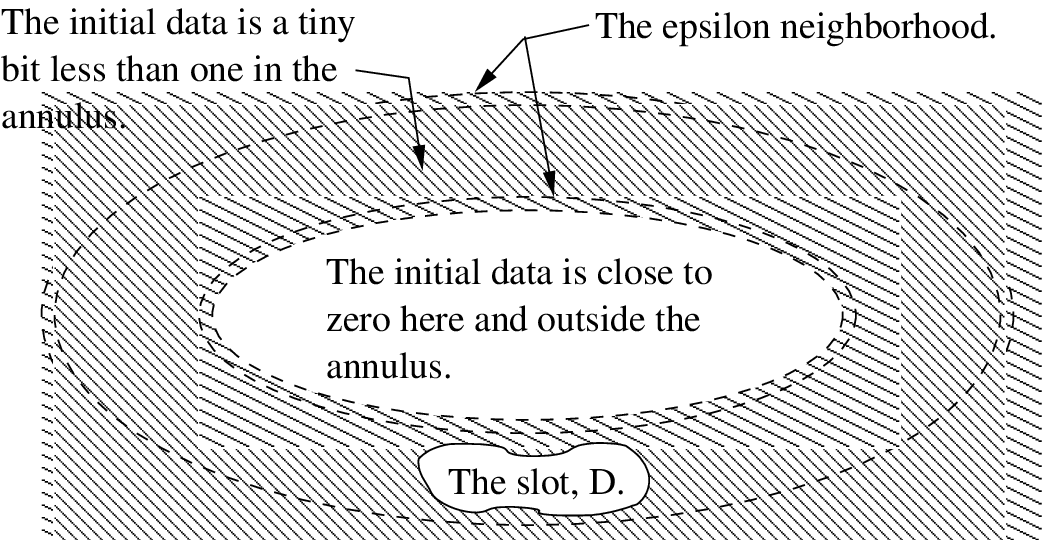}

\


\newsec{Continuity and dealing with all times}{CDWAT}
We wish now to improve the results from the previous section by extending
them from almost every time to every time.  Indeed, this prevents corners
from arising as ``transitions'' from the situations which are permissible
for sets of time with positive measure.  (Consider for example the sets
$S(t) := \{ (x,y) \; : \; xy < t - 1$ for the interval $t \in [0, 2].$  By
the results of the last section, it would still be possible for $S(t)$ to
be the diffusive region for a Hele-Shaw flow, since the corner only occurs
when $t$ belongs to the zero measure set $\{ 1 \}.$)  Of course in this 
entire section we make the standing assumption that $u_I$ is nondegenerate 
initial data (see Equation\refeqn{ndsit}\!).  We start with a simple lemma
summarizing some of the regularity we have for $W.$

\begin{lemma}[Regularity for $W$]   \label{RegW}
$W(x,t)$ is continuous in space and continuous and convex and nondecreasing 
in time.  For almost every time, $T,$ $W(x, T)$ satisfies 
Equation\refeqn{WobProb}with boundary data
\begin{equation}
    W(x, T) = p(x)T \ \ \ \ \ x \in \partial D \;.
\label{eq:Wbdrydata}
\end{equation}
\end{lemma}
\begin{pf}
By using the maximum principle together with the fact that the diffusive
region increases with time, we see that $V(x, t)$ must be an increasing
function of time.  Now convexity of $W$ follows from this fact and from
the definition of $W.$  The rest of the lemma follows immediately from 
the definition of $W$ and the spatial continuity of $V(x,t)$ along with
its boundary data on $\partial D.$
\newline Q.E.D. \newline
\end{pf}
\begin{theorem}[Measure of the diffusive region]   \label{medire}
There exists a constant $C = C(n, \alpha, \partial D)$ such that if 
$0 \leq t - s \leq 1,$ then
\begin{equation}
   |A(t) \setminus A(s)| \leq C ||p||_{C^{2,\alpha}(\partial D)} 
       \left( \frac{t - s}{1 - \lambda} \right) \;.
\label{eq:mcwrtt}
\end{equation}
(Recall that $A(t)$ is the diffusive region at time $t,$ and for
$S \subset \R^n,$ we let $|S|$ denote the Lebesgue n-dimensional measure
of $S.$)
\end{theorem}
\begin{pf}
We adapt the proof of Theorem 4.1 of \cite{Bl} to the current setting.
Fix $t$ and $s.$  Because the diffusive regions are nested, and because 
Equation\refeqn{WobProb}holds for almost every time, without loss of generality
we can assume that it holds for both $t$ and $s.$  Let 
$L := A(t) \setminus A(s),$ let $\Psi(x) := W(x,t) - W(x,s)$ and observe that
$\Delta \Psi(x) = \chisub{L} (1 - u_{I}) \geq 0,$ and 
\begin{equation}
\begin{array}{ll}
    \Psi(x) = 0 \ & x \in FB(t), \\
    \Psi(x) = (t - s)p(x) \ & x \in \partial D \;.
\end{array}
\label{eq:bvd4wsmwt}
\end{equation}
From this fact and by the weak maximum principle, it follows that
\begin{equation}
0 \leq W(x,t) - W(x,s) \leq (t - s)||p||_{L^{\infty}(\partial D)} 
\label{eq:beW}
\end{equation}
for all $x \in A(t).$
(Nonnegativity is actually a consequence of the previous lemma.)

Observe that $W(x,t) - W(x,s)$ is harmonic within $A(s)$ so that
\begin{alignat*}{1}
0 &= \int_{A(s)} \Delta(W(x,t) - W(x,s)) dx \\
  &= \int_{\partial D} \frac{\partial}{\partial \nu} (W(x,t) - W(x,s)) d\mathcal{H}^{n-1}
   - \int_{FB(s)} \frac{\partial}{\partial \nu} (W(x,t) - W(x,s)) d\mathcal{H}^{n-1} \\
  &= \int_{\partial D} \frac{\partial}{\partial \nu} (W(x,t) - W(x,s)) d\mathcal{H}^{n-1}
   - \int_{FB(s)} \frac{\partial}{\partial \nu} W(x,t) d\mathcal{H}^{n-1} 
\end{alignat*}
Now by using boundary regularity for harmonic functions combined with 
Equation\refeqn{beW}we can conclude
\begin{equation}
   \left| \int_{\partial D} \frac{\partial}{\partial \nu} 
          (W(x,t) - W(x,s)) d\mathcal{H}^{n-1} \right| \leq 
   C(n, \alpha, \partial D) (t - s)||p||_{C^{2,\alpha}(\partial D)}
\label{eq:ellregest}
\end{equation}
By combining this fact with the last computation, we conclude that
\begin{equation}
   \left| \int_{FB(s)} \frac{\partial}{\partial \nu} W(x,t) d\mathcal{H}^{n-1} \right|
   \leq C(n, \alpha, \partial D) (t - s) ||p||_{C^{2,\alpha}(\partial D)} \;.
\label{eq:bdryestI}
\end{equation}
On the other hand we have
\begin{alignat*}{1}
   (1 - \lambda)|L| 
      &\leq \int_{L} (1 - u_{I}) \; dx \\
      &= \int_{L} \Delta W(x,t) \; dx \\
      &= \int_{\partial L} \frac{\partial}{\partial \nu} W(x,t) d\mathcal{H}^{n-1} \\
      &= \int_{FB(s)} \frac{\partial}{\partial \nu} W(x,t) d\mathcal{H}^{n-1}
\end{alignat*}
which we can combine with Equation\refeqn{bdryestI}to give us what we need.
\newline Q.E.D. \newline
\end{pf}
\begin{corollary}[Continuity in $L^p$]  \label{CLp}
Under the assumptions made at the beginning of this section, the map from $t$ 
to the function
$$\chisub{\{ W(x,t) > 0\}} (1 - u_{I})$$
is a continuous function from $\R$ into $L^{p}(D^{c})$ for $1 \leq p < \infty.$
\end{corollary}
\begin{corollary}[Spatial regularity for every time]   \label{HOP}
All of the results of Theorems\refthm{RFBS}and\refthm{betreg}hold for every time.
\end{corollary}
\begin{pf}
It suffices to show that Equation\refeqn{WobProb}holds for all time.  Fix
$\tilde{t} > 0$ and let $t_n \rightarrow \tilde{t}$ with $t_n$ chosen so that
Equation\refeqn{WobProb}holds at each $t_n.$  Now take a ball, $B_R$ which is large
enough to contain $A(\tilde{t})$ in its interior, and let 
$\Omega := B_R \setminus D.$  For each $n$ we let $w_n(x)$ solve the boundary value
problem
\begin{equation}
   \begin{array}{ll}
            \Delta w_n(x) = \chisub{A(t_n)}(1 - u_{I}) & \ \text{in}
               \ \ \Omega \vspace{.05in} \\
            w_n(x) = t_n p(x) & \ \text{on} \ \ \partial D \\
            w_n(x) = 0 & \ \text{on} \ \ \partial B_R \;.
   \end{array}
\label{eq:bvpop}
\end{equation}
By standard uniqueness results, $w_n(x) \equiv W(x, t_n),$ as they satisfy the same
boundary value problem.  By standard elliptic regularity theory, since the boundary
data on $\partial D$ will converge to $\tilde{t} p(x)$ and by the last corollary the 
right hand side of the equation will converge in $L^{p}(\Omega)$ to 
$\chisub{A(\tilde{t})}(1 - u_{I}),$ we can conclude that $w_n$ will converge to a 
function $\tilde{w}$ which satisfies
\begin{equation}
   \begin{array}{ll}
            \Delta \tilde{w}(x) = \chisub{A(\tilde{t})}(1 - u_{I}) & \ \text{in}
               \ \ \Omega \vspace{.05in} \\
            \tilde{w}(x) = \tilde{t} p(x) & \ \text{on} \ \ \partial D \\
            \tilde{w}(x) = 0 & \ \text{on} \ \ \partial B_R \;.
   \end{array}
\label{eq:bvpoplim}
\end{equation}
On the other hand, $w_n(x) = W(x, t_n)$ converges to $W(x, \tilde{t})$ by the continuity
of $W$ in time.
\newline Q.E.D. \newline
\end{pf}


\newsec{Appendix}{App}
Here we collect some facts we need to construct our subsolutions in the proof of
Theorem\refthm{NMFB}\!.  We let $u(r; \alpha, \beta)$ denote the solution to:
\begin{equation}
  \begin{array}{rll}
    \Delta_x u(r; \alpha, \beta) \!
    &= 0 \ \ \
    & \text{in} \ B_{1 + \alpha + \beta} \setminus B_{\alpha} \\
    \ \\
    u(r; \alpha, \beta) \! &= 1 \ \ \ & \text{on} \ \partial B_{\alpha} \\
    \ \\
    u(r; \alpha, \beta) \! &= 0 \ \ \ & \text{on} \ \partial B_{1 + \alpha + \beta} \\
  \end{array}
\label{eq:TestU}
\end{equation}
and we let $v(r; \alpha, \beta)$ denote the solution to:
\begin{equation}
  \begin{array}{rll}
    \Delta v(r; \alpha, \beta) \!
    &= 2n \ \ \
    & \text{in} \ B_{1 + \alpha + \beta} \setminus B_{\alpha} \\
    \ \\
    v(r; \alpha, \beta) \!
    &= 0 \ \ \
    & \text{on} \ \partial \{ B_{1 + \alpha + \beta} \setminus B_{\alpha} \} \;.
  \end{array}
\label{eq:TestV}
\end{equation}
We will always assume that $\alpha \geq 1$ and $0 \leq \beta \leq 1.$

\begin{lemma}[Explicit Forms of Our Comparison Functions]  \label{RHFAR}
If $n > 2,$ then $u$ has the explicit form
\begin{equation}
   u(r; \alpha, \beta) = 
      \frac{r^{2 - n} - (1 + \alpha + \beta)^{2 - n}}
      {\alpha^{2 - n} - (1 + \alpha + \beta)^{2 - n}} \;,
\label{eq:Uform}
\end{equation}
and $v$ has the explicit form
\begin{equation}
   v(r; \alpha, \beta) = (r^2 - \alpha^2) + 
      \frac{(\alpha^2 - (1 + \alpha + \beta)^2)(r^{2 - n} - \alpha^{2 - n})}
           {(1 + \alpha + \beta)^{2 - n} - \alpha^{2 - n}} \;.
\label{eq:Vform}
\end{equation}
If $n = 2,$ then $u$ has the explicit form
\begin{equation}
   u(r; \alpha, \beta) = \left[ \log \left( \frac{r}{1 + \alpha + \beta} \right) \right] /
                    \left[ \log \left( \frac{\alpha}{1 + \alpha + \beta} \right) \right] \;,
\label{eq:Uform2}
\end{equation}
and $v$ has the explicit form
\begin{equation}
   v(r; \alpha, \beta) = (r^2 - \alpha^2) + 
   \frac{ \log \left( \displaystyle{ \frac{r}{\alpha} } \right) 
   (\alpha^2 - (1 + \alpha + \beta)^2)}{\log \left( \displaystyle{
      \frac{1 + \alpha + \beta}{\alpha} } \right) } \;.
\label{eq:Vform2}
\end{equation}
\end{lemma}

\begin{lemma}[Derivatives on the Outer Boundaries]  \label{DOB}
If $n > 2,$ then
\begin{equation}
   u_r(1 + \alpha + \beta; \alpha, \beta) = \frac{(2 - n)(1 + \alpha + \beta)^{1 - n}}
                     {\alpha^{2 - n} - (1 + \alpha + \beta)^{2 - n}} < 0 \;,
\label{eq:urfb}
\end{equation}
and
\begin{equation}
     v_r(1 + \alpha + \beta; \alpha, \beta) = 2(1 + \alpha + \beta) +
 (n - 2)(1 + \alpha + \beta) \frac{(1 + \alpha + \beta)^2 - \alpha^2}
                                  {(1 + \alpha + \beta)^{2 - n} - \alpha^{2 - n}} > 0 \;.
\label{eq:vrfb}
\end{equation}
If $n = 2,$ then
\begin{equation}
   u_r(1 + \alpha + \beta; \alpha, \beta) = \frac{1}
  {(1 + \alpha + \beta) \log \left( \displaystyle{\frac{\alpha}{1 + \alpha + \beta}} \right)} < 0 \;,
\label{eq:u2rfb}
\end{equation}
and
\begin{equation}
     v_r(1 + \alpha + \beta; \alpha, \beta) = 2(1 + \alpha + \beta) +
        \frac{\alpha^2 - (1 + \alpha + \beta)^2}{(1 + \alpha + \beta) \log \left(
        \displaystyle{\frac{1 + \alpha + \beta}{\alpha}} \right)}
                                   > 0 \;.
\label{eq:v2rfb}
\end{equation}
Also, for any $n \geq 2,$ we have
\begin{equation}
   \lim_{\alpha \rightarrow \infty} u_r(1 + \alpha + \beta; \alpha, \beta) =
   \frac{-1}{1 + \beta} < 0 \;,
\label{eq:urfbLIM}
\end{equation}
and
\begin{equation}
   \lim_{\alpha \rightarrow \infty} v_r(1 + \alpha + \beta; \alpha, \beta) =
   n(1 + \beta) > 0 \;.
\label{eq:vrfbLIM}
\end{equation}
\end{lemma}
\begin{corollary}[Bounds on the Outer Boundaries]  \label{BOB}
For $0 \leq \beta \leq 1$ and $\alpha \geq 1,$ there exist constants $\gamma_i$
which are all independent of $\alpha$ and $\beta$ such that
\begin{equation}
   -\gamma_1 \leq u_r(1 + \alpha + \beta; \alpha, \beta) \leq -\gamma_2 < 0 \;,
\label{eq:uHG}
\end{equation}
and
\begin{equation}
   0 < \gamma_3 \leq v_r(1 + \alpha + \beta; \alpha, \beta) \leq \gamma_4 < \infty \;.
\label{eq:vHG}
\end{equation}
\end{corollary}
\begin{pf}
The proof for $u_r$ is essentially the same as the proof for $v_r,$ so we will only
deal with $u_r.$  Because the limit as $\alpha \rightarrow \infty$ is strictly negative
by the previous lemma there is a large $\alpha_0$ such that for $\alpha \geq \alpha_0$
we have the desired lower and upper bounds.  Next we simply use compactness for the 
rest of the strip.
\newline Q.E.D. \newline
\end{pf}

\end{document}